# A COMPLEMENT TO LE CAM'S THEOREM

By Mark G. Low[1] and Harrison H. Zhou

*University of Pennsylvania and Yale University*

This paper examines asymptotic equivalence in the sense of Le Cam between density estimation experiments and the accompanying Poisson experiments. The significance of asymptotic equivalence is that all asymptotically optimal statistical procedures can be carried over from one experiment to the other. The equivalence given here is established under a weak assumption on the parameter space $\mathcal{F}$. In particular, a sharp Besov smoothness condition is given on $\mathcal{F}$ which is sufficient for Poissonization, namely, if $\mathcal{F}$ is in a Besov ball $B_{p,q}^{\alpha}(M)$ with $\alpha p > 1/2$. Examples show Poissonization is not possible whenever $\alpha p < 1/2$. In addition, asymptotic equivalence of the density estimation model and the accompanying Poisson experiment is established for all compact subsets of $C([0,1]^m)$, a condition which includes all Hölder balls with smoothness $\alpha > 0$.

**1. Introduction.** A family of probability measures $E = \{P_\theta : \theta \in \Theta\}$ defined on the same $\sigma$-field is called a statistical model. Le Cam [8] defined a distance $\Delta(E, F, \Theta)$ between $E$ and another model $F = \{Q_\theta : \theta \in \Theta\}$ with the same parameter set $\Theta$. For bounded loss functions, if $\Delta(E, F, \Theta)$ is small, then to every statistical procedure for $E$ there is a corresponding procedure for $F$ with almost the same risk function.

Le Cam [9] used the deficiency distance between the experiment $E_n$ with $n$ i.i.d. observations and the experiment $E_{n+r}$ with $n + r$ i.i.d. observations as a measure of the amount of information in the additional observations. It was shown that the deficiency distance $\Delta(E_n, E_{n+r}, \mathcal{F})$ can be bounded by

$$\Delta(E_n, E_{n+r}, \mathcal{F}) \leq \sqrt{8r\beta_n}, \tag{1}$$

where $\beta_n$ is the minimax risk for estimating a density function under squared Hellinger distance based on the experiment $E_n$. For any two measures $P$

Received July 2005; revised June 2006.
[1]Supported in part by NSF Grant DMS-03-06576.

and $Q$ the Hellinger distance $H(P,Q)$ is defined by $H^2(P,Q) = \int(\sqrt{dP} - \sqrt{dQ})^2$. For regular parametric models $\beta_n$ is of order $n^{-1}$ and Le Cam's upper bound for $\Delta(E_n, E_{n+r}, \mathcal{F})$ is then $C(r/n)^{1/2}$ for some $C > 0$. This bound was furthered improved in Mammen [11] to $Cr/n$ once again for the case of regular parametric models.

As pointed out by Le Cam [9] the information content in additional observations is connected to the "technical device which consists in replacing the fixed sample size $n$ by a Poisson variable $N$." More specifically throughout this paper we shall consider the following two experiments.

*Density estimation experiment.*

(2) $\qquad E_n : y_1, y_2, \ldots, y_n$ i.i.d. with density $f$.

*Poisson experiment.*

(3) $\qquad F_n : \underline{x}_n(\cdot)$, a Poisson process with intensity measure $nf$.

Equivalently, the Poisson experiment corresponds to observing a Poisson random variable $N$ with expectation $n$ and then independently of $N$ observing $y_1, y_2, \ldots, y_N$ i.i.d. with density $f$. For both experiments $f$ is an unknown density and $f \in \mathcal{F}$ a given parameter space and we shall say that Poissonization is possible if $\Delta(E_n, F_n, \mathcal{F}) \to 0$. Le Cam [9] showed $\Delta(E_n, F_n, \mathcal{F}) \leq Cn^{-1/4}$ for regular parametric models and he also gave the following general result.

PROPOSITION 1 (Le Cam). *Suppose that there is a sequence of estimators $\hat{f}_n$ based on either the density estimation model $E_n$ or the Poisson process model $F_n$ satisfying*

(4) $$\sup_{f \in \mathcal{F}} E_f n^{1/2} H^2(\hat{f}_n, f) \to 0.$$

*Then*

$$\Delta(E_n, F_n, \mathcal{F}) \to 0.$$

It should be noted that the condition (4) is quite a strong assumption. However for Hölder spaces defined on the unit interval with $\alpha > 1/2$, Yang and Barron [15] showed that there is an estimator for which $\beta_n^* = o(1/\sqrt{n})$ and in this case it follows from Proposition 1 that $\Delta(E_n, F_n, \mathcal{F}) \to 0$. We should also note that Poissonization is not always possible. Le Cam [9] does give an example of a parameter space for which $\Delta(E_n, F_n, \mathcal{F}) \not\to 0$. However the parameter space used for this counterexample is so "large" that there does not even exist a uniformly Hellinger consistent estimator over this parameter space. There is thus a considerable gap in the condition given in Proposition 1 which guarantees that Poissonization is possible and this example for which Poissonization fails. The present paper aims to at least partially fill this gap.

In the last decade much progress has been made in bounding the deficiency distance between nonparametric models and Gaussian white noise models. In particular, theory has been developed for nonparametric density estimation models in Nussbaum [13], nonparametric regression in Brown and Low [2], generalized linear models in Grama and Nussbaum [5], for nonparametric autoregression in Milstein and Nussbaum [12] and for spectral density models in Golubev, Nussbaum and Zhou [4]. In all of this work asymptotic equivalence is established under particular smoothness assumptions on a nonparametric model which in terms of a Hölder smoothness condition for functions defined on $[0,1]$ corresponds to an assumption of at least $\alpha > 1/2$. As noted earlier the condition $\alpha > 1/2$ is exactly the minimal Hölder smoothness for which the assumption in Proposition 1 holds. Moreover, for the cases of nonparametric regression and nonparametric density estimation these models have been shown to be not equivalent to the corresponding white noise with drift model when $\alpha < 1/2$. See Brown and Low [2] and Brown and Zhang [3]. A corresponding theory has not yet been developed which explains when Poissonization is possible for such nonparametric models.

The focus of the present paper is to develop such a theory. We start, in Section 2, by giving some further examples where Poissonization is not possible. These examples are interesting because the parameter spaces used in these examples are much smaller than the one given in Le Cam [9]. In particular, the minimax rates of convergence under squared Hellinger distance can be of order $n^{-\gamma}$ with $\gamma$ arbitrarily close to $1/2$. Thus in terms of Hellinger distance the sufficient condition given in Proposition 1 cannot be improved. However, examples of parameter spaces are also given for which the minimax Hellinger distance converges to zero at a rate $n^{-\gamma}$ with $\gamma$ arbitrarily close to zero but where Poissonization holds. Taken together these results show that Hellinger distance cannot fully explain when Poissonization is possible.

The focus of Section 3 is on developing an alternative sufficient condition which guarantees that Poissonization is possible. A sequence of loss functions is introduced which are bounded between a chi-square distance and a squared Hellinger distance. It is shown that if there exists a sequence of uniformly consistent estimators under this sequence of loss functions then Poissonization is possible and $\Delta(E_n, E_{n+D\sqrt{n}}, \mathcal{F}) \to 0$ for every $D > 0$ (see Theorem 2). In particular, in contrast to the theory for Gaussian equivalence Poissonization is possible over all Hölder balls with $\alpha > 0$.

The theory also allows for a characterization of the Besov spaces for which Poissonization is possible. Under the sequence of losses defined in Section 3, a uniformly consistent sequence of estimators is constructed for Besov spaces with parameters $\alpha p > 1/2$, demonstrating that in these cases Poissonization is possible. On the other hand, the examples given in Section 2 show that Poissonization is not possible for Besov spaces with $\alpha p < 1/2$.

**2. Examples where Poissonization is not possible.** As mentioned in the introduction Le Cam [9] gave an example of a parameter set and a statistical problem which showed that the deficiency distance between i.i.d. observations and the Poissonized version of this experiment does not go to zero. In this example the observations have support on the unit interval but the parameter space, say $\mathcal{F}$, is not precompact in Hellinger distance. In fact, in his example uniformly consistent estimators under Hellinger distance do not even exist. Every sequence of estimators $\hat{f}_n$ satisfies

$$\liminf_{n \to \infty} \sup_{f \in \mathcal{F}} EH^2(\hat{f}_n, f) > 0.$$

This is the only example in the literature that we are aware of for which $\Delta(E_n, F_n, \mathcal{F}) \nrightarrow 0$.

In this section it is shown that there are "much smaller" parameter spaces for which Poissonization is not possible. In each of these examples the parameter space is "much smaller" than that given by Le Cam and in particular is compact under Hellinger distance. Moreover, for every $r < \frac{1}{2}$ an example is given of a parameter space $\mathcal{F}$ for which the minimax risk satisfies

(5) $$\inf_{\hat{f}_n} \sup_{f \in \mathcal{F}} E n^r H^2(\hat{f}_n, f) \to 0$$

but where Poissonization is not possible.

2.1. *Besov spaces.* The counterexamples we provide in this section are given for Besov spaces. These spaces occupy a central role in much of the recent nonparametric function estimation literature. Besov spaces also arise naturally in equivalence theory. In Brown, Carter, Low and Zhang [1] they were used to characterize when the density estimation model is asymptotically equivalent to a Gaussian process model.

Let $J_{j,k}$ be the averaging operator

$$J_{j,k}(f) = k \int_{(j-1)/k}^{j/k} f(x)\, dx$$

and define the piecewise constant approximation $\bar{f}_{(k)}$ by

(6) $$\bar{f}_{(k)} = \sum_{j=1}^{k} J_{j,k}(f) 1_{[(j-1)/k, j/k)}.$$

Then for each function $f$ on $[0,1]$ the Besov norm is given by

$$\|f\|_{\alpha,p,q} = \left\{\left|\int f(x)\,dx\right|^q + \sum_{i=0}^{\infty}(2^{i\alpha}\|\bar{f}_{(2^{i+1})} - \bar{f}_{(2^i)}\|_p)^q\right\}^{1/q}, \tag{7}$$

and the Besov balls can then be defined by

$$B_{p,q}^{\alpha}(M) = \{\|f\|_{\alpha,p,q} \leq M\}.$$

Under squared Hellinger distance rate optimal estimators have been constructed in Yang and Barron [15]. It was shown that

$$C_1(n\log n)^{-2\alpha/(2\alpha+1)} \leq \inf_{\hat{f}} \sup_{f\in\mathcal{F}} EH^2(\hat{f},f)$$
$$\leq C_2 n^{-2\alpha/(2\alpha+1)}(\log n)^{1/(2\alpha+1)}, \tag{8}$$

when $\alpha + 1/2 - 1/p > 0$, $p \geq 1$ and $q \geq 1$. We should note that it immediately follows from the sufficient condition of Le Cam that $\Delta(E_n, F_n, \mathcal{F}) \to 0$ when the parameter space $\mathcal{F} = B_{p,q}^{\alpha}(M)$ and $\alpha > \frac{1}{2}$. In the counterexamples that follow $\alpha p < 1/2$. In Section 4 it is shown that Poissonization is possible whenever $\alpha p > 1/2$.

2.2. *Counterexamples.* In order to show that $\Delta(E_n, F_n, \mathcal{F}) \not\to 0$ it suffices to exhibit a sequence of statistical problems with bounded loss functions such that the Bayes risks are asymptotically different. This approach taken by Le Cam [9] and Brown and Zhang [3] requires the specification of a sequence of decision problems along with a particular sequence of priors for which the Bayes risks are asymptotically different.

We adopt the same general strategy. First we shall provide a description of the sequence of priors and then we shall turn to the particular decision problems. The priors we shall use correspond to uniform priors placed on only a finite number of functions in $\mathcal{F}$. For this reason it is convenient to specify these priors by first describing the set of points on which the priors are supported. For $n \geq 1$, let $I_n$ be the collection of intervals $[\frac{i-1}{n}, \frac{i}{n})$, $i = 1, 2, \ldots, n$. Now for $1/2 < \beta < 1$ define

$$\mathcal{F}_{\beta,n} = \left\{f:[0,1]\to R : f = 0 \text{ on } n^\beta \text{ intervals in } I_n \text{ and } f = \frac{n}{n-n^\beta} \text{ otherwise}\right\}.$$

It is simple to check by computing the Besov norms that for any $M > 1$, $\mathcal{F}_{\beta,n} \subset B_{p,q}^{\alpha}(M)$ for $n$ sufficiently large, whenever $\alpha p < 1 - \beta$. The priors that we shall use correspond to uniform priors on these sets where the $\beta$ is chosen so that $\alpha p < 1 - \beta$.

We now turn to a description of a collection of decision problems. For a given known $m$ which may depend on $n$, using either i.i.d. data or the Poissonized data, we wish to name exactly $m$ intervals of length $1/n$ where the function is not zero. More specifically, we must list $m$ intervals of the form $[(i-1)/n, i/n)$, where the integer $i$ satisfies $1 \leq i \leq n$. For this problem we impose the following loss function. If we name $m$ such intervals correctly then the loss is zero. If we even make one mistake the loss is 1. The difficulty of this problem depends strongly on the magnitude of $m$ as well as on the value of $\beta$. For example, if $m$ is small then just random guessing of such intervals usually results in zero loss since the function takes on a nonzero value on most of the intervals. The problem becomes difficult when $m$ is large. This idea can be developed further as follows. Let $K_E$ and $K_F$ be equal to the number of intervals containing at least one observation based on the density estimation model and Poisson process model, respectively. Then it is easy to calculate the expectations and variances of these random variables. Taylor series expansions yield

$$E(K_E) = n(1 - e^{-1}) + n^\beta(-1 + 2e^{-1}) + O\left(\frac{n^{2\beta}}{n}\right)$$

and

$$E(K_F) = n(1 - e^{-1}) + n^\beta(-1 + 2e^{-1}) + O\left(\frac{n^{2\beta}}{n}\right).$$

Likewise the variances satisfy

$$\mathrm{Var}(K_E) = n(e^{-1} + o(1)),$$
$$\mathrm{Var}(K_F) = n((1 - e^{-1})e^{-1} + o(1)).$$

Our counterexamples are constructed by choosing a value of $m$ where the variability of $K_E$ and $K_F$ plays an important role in the difficulty of the problem. That is, we shall take $m$ to be equal to the expected value of $K_E$ minus a small multiple of the standard deviation of $K_E$ or $K_F$. Specifically, set $m = n(1 - e^{-1}) + n^\beta(-1 + 2e^{-1}) - \sqrt{n}$. For such an $m$ the chance that $K_E < m$ differs significantly from the chance that $K_F < m$.

It is convenient to recast the problem in the following way. Note that once we have decided on a set of $m$ intervals where the function is not equal to zero we can take the subset of $\mathcal{F}_{\beta,n}$ which contains all such functions with this property. Call this set of functions $S$. The loss associated to $S$ is then

$$L(f, S) = \begin{cases} 0, & f \in S, \\ 1, & \text{otherwise.} \end{cases}$$

Recast in this manner the problem is thus to select $S$, a subset of $\mathcal{F}_{\beta,n}$, each satisfying $f = n/(n - n^\beta)$ on those $m$ intervals. That is, the set $S$ is equal to the collection of functions in $\mathcal{F}_{\beta,n}$ which take on the value $f = n/(n - n^\beta)$ on the $m$ intervals.

As mentioned at the beginning of this section, in order to show that the i.i.d. observations and the Poissonized version are asymptotically nonequivalent we shall show that the Bayes risks for these two problems are different when we put a uniform prior on $\mathcal{F}_{\beta,n}$. A Bayes solution to each of these problems is straightforward. In $E_n$ when $K_E \geq m$ or in $F_n$ when $K_F \geq m$ the selection of $S$ is easy. In these cases we know $m$ intervals where the function is equal to $f = n/(n - n^\beta)$ and we can just take $S$ to be a set of functions in $\mathcal{F}_{\beta,n}$ with this property. The loss suffered in this case is clearly 0. If $K_E < m$ or $K_F < m$ we need to choose an additional $m - K_E$ or $m - K_F$ intervals in order to construct $S$. A Bayes rule for doing this is to select $m - K_E$ or $m - K_F$ additional intervals randomly from the remaining $n - K_E$ or $n - K_F$ intervals based on the uniform prior over these intervals. Writing $K$ for either $K_E$ or $K_F$ we see that the expected loss for these problems given the value of $K$ when $K < m$ is just 1 minus the chance that, when picking $m - K$ balls out of an urn with $n - n^\beta - K$ black balls and a total of $n - K$ balls, each ball chosen is black. The chance that this occurs is just

$$1 - \binom{n - n^\beta - K}{m - K} \bigg/ \binom{n - K}{m - K}.$$

Hence the Bayes risk for these problems can be written as

$$R_n = E\left[1 - \binom{n - n^\beta - K}{m - K} \bigg/ \binom{n - K}{m - K}\right] I_{\{K < m\}},$$

where $K$ is either $K_E$ or $K_F$.

In the Appendix we prove the following lemma.

LEMMA 1. *With $R_n$ defined above in both the density estimation setting $E_n$ and the Poisson process setting $F_n$,*

(9) $$R_n = P(K < m) + o(1),$$

*where $K = K_E$ for $E_n$, or $K_F$ for $F_n$.*

It is then easy to see that the value of $R_n$ is asymptotically different for $E_n$ and $F_n$. For $E_n$ note that the central limit theorem (CLT) for the occupancy problem (see Kolchin, Sevast'yanov and Chistyakov [7]) shows that

$$P(K_E < m) \to \Phi(-\sqrt{e}),$$

and the usual CLT yields

$$P(K_F < M) \to \Phi(-\sqrt{e}/\sqrt{1 - 1/e}),$$

where $\Phi$ is the cumulative distribution function of the standard normal distribution. Hence the Bayes risks for this problem differ asymptotically.

Now consider the Besov space $B_{p,q}^\alpha(M)$ with $M > 1$ and $\alpha p < 1/2$. Then take $1/2 < \beta < 1 - \alpha p$. It then follows that for sufficiently large $n$, $\mathcal{F}_{\beta,n} \subset B_{p,q}^\alpha(M)$. Since, as we have just shown, there is a sequence of priors on $\mathcal{F}_{\beta,n}$ which have different asymptotic Bayes risks for i.i.d. data and Poisson data, the same is trivially true for $B_{p,q}^\alpha(M)$. This in turn shows that the deficiency distance does not tend to zero. The consequence of these results for asymptotic equivalence can then be summarized in the following theorem.

THEOREM 1. *Suppose $\alpha p < \frac{1}{2}$ and $M > 1$. Then*

$$\Delta(E_n, F_n, B_{p,q}^\alpha(M)) \nrightarrow 0.$$

REMARK 1. Note that choosing $p = 1$ and some $\alpha < 1/2$, the results of Yang and Barron [15] given in (8) show that for any algebraic rate of convergence slower than $n^{-1/2}$ there are Besov parameter spaces with at least this rate of convergence, under squared Hellinger distance, where Poissonization is not possible.

**3. Asymptotic equivalence under a general assumption.** In the previous section examples were presented where the rate of convergence under squared Hellinger distance is arbitrarily close to $1/2$ but where Poissonization is not possible. It follows that any weakening of Le Cam's assumptions for Poissonization must involve something other than Hellinger distance.

In this section it is shown that Poissonization is possible under a condition which substantially improves on the sufficient condition given in Le Cam [9]. In particular, for all Hölder balls on the unit interval with arbitrary smoothness $\alpha > 0$ Poissonization is possible although the sufficient condition of Le Cam given in Proposition 1 shows Poissonization is possible only if $\alpha > 1/2$.

Considerable insight into a comparison of the two experiments $E_n$ and $F_n$ can be gained by the following simple observation which also greatly simplifies the analysis. Consider two Poisson experiments, $F_{n-m}$ and $F_{n+m}$ where $m = cn^\gamma$. If $\gamma > \frac{1}{2}$, then with probability approaching one the number of observations from $F_{n-m}$ is less than $n$ whereas the number of observations from $F_{n+m}$ is larger than $n$. It is then easy to check that asymptotically $E_n$ is at least as informative as $F_{n-m}$ and $F_{n+m}$ is at least as informative as $E_n$. In fact, by taking $\gamma = \frac{1}{2}$ and $c$ sufficiently large simple bounds based on the chance that there are more or less than $n$ observations lead to bounds on how much more or less informative these experiments can be. The following lemma captures these ideas. The proof can be found in the Appendix.

LEMMA 2. *Suppose that for each $D > 0$*

$$\lim_{n \to \infty} \Delta(F_n, F_{n+D\sqrt{n}}, \mathcal{F}) = 0.$$

*Then*

$$\lim_{n \to 0} \Delta(E_n, F_n, \mathcal{F}) \to 0.$$

Thus in order to show that $\Delta(E_n, F_n, \mathcal{F}) \to 0$ we can focus on measuring the deficiency distance between two Poisson process experiments. It should be noted that some insight into the deficiency distance between two Poisson experiments is provided by the following general bound given in Le Cam [10]:

$$\Delta(F_n, F_{n+m}, \mathcal{F}) \leq \frac{m}{\sqrt{2n}};$$

but clearly this bound does not suffice in the present context.

It is useful to recall the Hellinger distance between any two Poisson processes with intensity $g$ and $h$. Write $P_f$ for the distribution of a Poisson process with intensity $f$. Then it follows from Le Cam ([10], page 160) that

$$H^2(P_g, P_h) = 2\left(1 - \exp\left(-\tfrac{1}{2}\int(\sqrt{g} - \sqrt{h})^2\right)\right)$$
$$\leq \int(\sqrt{g} - \sqrt{h})^2.$$

In particular, the following upper bound holds for the Hellinger distance between Poisson processes with intensities $nf$ and $(n+m)f$:

$$H^2(P_{nf}, P_{(n+m)f}) \leq \int(\sqrt{nf} - \sqrt{(n+m)f})^2 = (\sqrt{n} - \sqrt{n+m})^2 \asymp \frac{m^2}{n}.$$

For this reason, to show $\Delta(F_n, F_{n+m}, \mathcal{F}) \to 0$ a randomization of the Poisson process with intensity $nf$ must be given which more closely matches that of the Poisson process with intensity $(n+m)f$. If we know that $f$ is in a neighborhood of a particular $f_0 \in \mathcal{F}$, this is easily accomplished by a superposition of the Poisson processes with intensities $nf$ and $mf_0$. For this new Poisson process we can calculate the Hellinger distance to yield

$$H^2(P_{nf+mf_0}, P_{(n+m)f}) \leq \frac{m^2}{n+m}\int \frac{(f - f_0)^2}{f + mf_0/(n+m)}.$$

In particular, if $m = D\sqrt{n}$ with $D > 1$ it immediately follows that

$$H^2(P_{nf+mf_0}, P_{(n+m)f}) \leq 2D^2 \int \frac{(f - f_0)^2}{f + n^{-1/2}f_0}.$$

The following result immediately follows.

LEMMA 3. *Set*

$$\mathcal{F}(f_0, c_n) = \left\{ f : \int \frac{(f - f_0)^2}{f + n^{-1/2} f_0} \leq c_n \right\}. \tag{10}$$

*Then*

$$\Delta(F_n, F_{n+D\sqrt{n}}, \mathcal{F}(f_0, c_n)) \leq 2D\sqrt{c_n}.$$

This lemma yields a general approach to giving a sufficient condition under which Poissonization is possible. Le Cam showed that in order to establish asymptotic equivalence for the whole parameter space it suffices to establish local asymptotic equivalence, as in Lemma 3, along with the existence of estimators which with probability tending to 1 localize you within such a neighborhood. In the present context it is natural to link the local parameter space around a given $f_0$ with the following loss function which also depends on $n$.

Let the loss $L_n$ be defined by

$$L_n(f, g) = \int \frac{(g - f)^2}{f + n^{-1/2} g} \, d\mu.$$

The following theorem, the proof of which is given in the Appendix, then gives a sufficient condition for Poissonization. This step is often called globalization.

THEOREM 2. *Let $\mathcal{F}$ be a parameter space that is separable under squared Hellinger distance. Fix $\varepsilon > 0$ and let $\hat{f}_{n,\varepsilon}$ be an estimator based on the model $E_n$. Suppose that $\hat{f}_{n,\varepsilon}$ satisfies*

$$\sup_{f \in \mathcal{F}} P_f \{ L_n(f, \hat{f}_{n,\varepsilon}) > \varepsilon \} \to 0. \tag{11}$$

*Then*

$$\Delta(E_n, F_n, \mathcal{F}) \to 0.$$

*In addition, we have $\Delta(E_n, E_{n+D\sqrt{n}}, \mathcal{F}) \to 0$ for every $D > 0$.*

Of course this theorem would not be particularly useful unless we were able to give some interesting examples under which (11) holds. However, before giving such examples it is worthwhile to note that although the loss function $L_n$ is not standard, it can be connected with squared Hellinger distance by using the inequalities

$$1 \leq \frac{(\sqrt{g} + \sqrt{f})^2}{f + n^{-1/2} g} \leq (1 + \sqrt{n}). \tag{12}$$

It then follows that

$$H^2(f,g) \leq L_n(f,g) = \int (\sqrt{g} - \sqrt{f})^2 \frac{(\sqrt{g}+\sqrt{f})^2}{f+n^{-1/2}g} \, d\mu \leq (1+\sqrt{n})H^2(f,g).$$

It thus immediately follows that any sequence of estimators which satisfies (4) also satisfies (11). Hence convergence under $L_n$ is weaker than Le Cam's condition.

REMARK 2. Let $\mathcal{F}$ be a compact subset of $C([0,1]^m)$ where $C([0,1]^m)$ is the collection of all continuous functions on the unit hypercube in $R^m$ with the $L_\infty$ norm as the measure of distance between functions. Standard arguments such as those found in Woodroofe [14] show that there exist estimators $\hat{f}_n$ such that for every $\varepsilon > 0$

$$\sup_{f \in \mathcal{F}} P(\|\hat{f}_n - f\|_\infty \geq \varepsilon) \to 0.$$

Define $\tilde{f}_n = \hat{f}_n 1(\hat{f}_n \geq 2\varepsilon)$. On the event $A_n = \{\|\hat{f}_n - f\|_\infty < \varepsilon\}$, we have $f(x) \leq \hat{f}_n(x) + |\hat{f}_n(x) - f(x)| \leq 3\varepsilon$ when $\tilde{f}_n = 0$ and $f(x) \geq \hat{f}_n(x) - |\hat{f}_n(x) - f(x)| \geq \varepsilon$ when $\hat{f}_n \geq 2\varepsilon$. It then follows that when $A_n$ occurs

$$L_n(f, \tilde{f}_n) \leq \int_{\{\hat{f}_n < 2\varepsilon\}} f + \int_{\{\hat{f}_n \geq 2\varepsilon\}} \frac{(\hat{f}_n - f)^2}{f}$$

$$\leq 3\varepsilon + \frac{\varepsilon^2}{\varepsilon} = 4\varepsilon.$$

Thus as $n \to \infty$

$$\sup_{f \in \mathcal{F}} P\{L_n(f, \tilde{f}_n) > 4\varepsilon\} \leq \sup_{f \in \mathcal{F}} P(A_n^c) \to 0.$$

It thus follows that Poissonization is possible in such cases. In particular, an example of compact subsets of $C([0,1]^m)$ is the set of Hölder balls $\mathcal{F} = \{f : |f(y) - f(x)| \leq M\|x - y\|^\alpha\}$ where $\|\cdot\|$ is the usual Euclidean norm on $R^m$.

REMARK 3. Suppose $\mathcal{F}$ is a compact subset of the space of functions on $\Omega \subset R^m$ under $L_2$ distance and that there is a $c > 0$ such that $f(x) \geq c$ for all $x \in \Omega$ and all $f \in \mathcal{F}$. Under this assumption, for any $\varepsilon > 0$ there is an estimator such that

$$\lim_{n \to \infty} \sup_f P_f(\|\hat{f}_n - f\|_2^2 > \varepsilon) = 0.$$

Note that

$$L_n(f,g) \leq \int \frac{(g-f)^2}{f} \leq \frac{1}{c}\|f - g\|_2^2$$

for any $g$ and $f$ in $\mathcal{F}$.

Thus for every $\varepsilon > 0$ we have

$$\sup_{f \in \mathcal{F}} P_f\{L_n(f, \hat{f}_n) > \varepsilon\} \leq \sup_{f \in \mathcal{F}} P_f(\|\hat{f}_n - f\|_2^2 > c\varepsilon) \to 0.$$

It follows that Poissonization is possible in such cases. In particular, any subset of Besov balls on the unit interval with $\alpha + 1/2 - 1/p > 0$ which have functions uniformly bounded away from 0 and above satisfies such a condition.

**4. Asymptotic equivalence for the unit interval.** In the previous section it was shown that the existence of consistent estimators under the loss function $L_n$ is sufficient for Poissonization and some examples of parameter spaces were given where such consistent estimators exist. In this section attention is focused on functions defined on the unit interval. Sufficient conditions on the parameter space $\mathcal{F}$ are given which guarantee the existence of uniformly consistent estimators under $L_n$, which in turn guarantees that Poissonization is possible.

Since the loss $L_n$ imposes a large penalty when the underlying function is close to zero but the estimator is not, it is natural to construct procedures which take on the value zero whenever it is suspected that the true function is close to zero. At the same time consistency under $L_n$ also requires the procedure to be close to the unknown function over most of the unit interval. This motivates the following simple modification of a histogram estimator for functions defined on the unit interval. We focus on the Poisson model. First consider the histogram estimator

$$\widehat{\bar{f}}_n(x) = \frac{k}{n} N_j, \qquad x \in [(j-1)/k, j/k), j = 1, 2, \ldots, k, \tag{13}$$

where $k = [n/\log^4 n]$ and where $N_j$ is the number of observations on the interval $[(j-1)/k, j/k)$. Note that $\widehat{\bar{f}}_n(x)$ defines a histogram on a very fine grid and that $\widehat{\bar{f}}_n(x)$ is an unbiased estimator of $\bar{f}_{(k)}$. The following modification of $\widehat{\bar{f}}_n(x)$ leads to a sequence of estimators which is often consistent under $L_n$.

Set $c_n = 1/\sqrt{\log n}$ and let $\hat{f}_n$ be defined by

$$\hat{f}_n(x) = \begin{cases} 0, & \widehat{\bar{f}}_n(x) < 2c_n, \\ 1/c_n, & \widehat{\bar{f}}_n(x) > 1/c_n, \\ \widehat{\bar{f}}_n(x), & \text{otherwise.} \end{cases} \tag{14}$$

The following theorem gives a structural condition on the parameter space $\mathcal{F}$ which guarantees that $\hat{f}_n$ is consistent under $L_n$.

THEOREM 3. *Let $\mathcal{F}$ be a collection of densities on the unit interval such that for some fixed $C > 0$, $\int f^2 \leq C$. Moreover suppose that*

$$\sup_{f \in \mathcal{F}} \mu\{x : |f(x) - \bar{f}(x)_{(k)}| > 1/\sqrt{\log k}\} = O(k^{-\delta}),$$

(15)

*for some $\delta > 1/2$.*

*Then for the Poisson model ([3](#)), the estimator $\hat{f}_n$ satisfies*

(16) $$\sup_{f \in \mathcal{F}} E_f L_n(f, \hat{f}_n) = o(1)$$

*and hence*

(17) $$\Delta(E_n, F_n, \mathcal{F}) \to 0.$$

*In particular ([16](#)) holds for Besov spaces $B_{p,q}^\alpha(M)$ whenever $\alpha p > 1/2$, $p \geq 1$ and $q \geq 1$.*

## APPENDIX

**A.1. Review of deficiency distance.** For any two experiments $E$ and $F$ with a common parameter space $\mathcal{F}$ the deficiency distance $\Delta(E, F, \mathcal{F})$ is defined by

$$\Delta(E, F, \mathcal{F}) = \max(\delta(E, F, \mathcal{F}), \delta(F, E, \mathcal{F})),$$

where

$$\delta(E, F, \mathcal{F}) = \inf_K \sup_{f \in \mathcal{F}} \|KP_f - Q_f\|_{\text{TV}},$$

where $K$ is a transition which is usually given by a Markov kernel. The triangle inequality

(18) $$\delta(E, G, \mathcal{F}) \leq \delta(E, F, \mathcal{F}) + \delta(F, G, \mathcal{F})$$

is used below in the proof of Lemma [2](#). Bounds between Hellinger distance and total variation immediately yield the following bound for the deficiency distance between the experiments $E = \{P_f, f \in \mathcal{F}\}$ and $F = \{Q_f, f \in \mathcal{F}\}$:

$$\Delta(E, F, \mathcal{F}) \leq 2 \inf_K \sup_{f \in \mathcal{F}} H^2(KP_f, Q_f).$$

**A.2. Proof of Lemma [1](#).** A proof is given only for the Poisson process model, as the proof for the density estimation model is similar. We have

$$R_n = E\left(1 - \frac{n - n^\beta - K_F}{n - K_F} \cdots \cdots \frac{n - n^\beta - K_F - (m - K_F - 1)}{n - K_F - (m - K_F - 1)}\right) I_{\{K_F < m\}}.$$

The general central limit theorem yields

$$\frac{K_F - EK_F}{\sqrt{\text{Var}(K_F)}} \rightsquigarrow N(0,1).$$

For any $0 < \varepsilon$ there are constants $k_1$ and $k_2$ such that

$$P(K_F < m - k_2\sqrt{n}) + P(m - k_1\sqrt{n} < K_F < m) \leq \varepsilon$$

for sufficiently large $n$. Simple calculation shows

$$\frac{n - n^\beta - K}{n - K} \cdots \frac{n - n^\beta - K - (m - K - 1)}{n - K - (m - K - 1)} \to 0$$

for $\beta > 1/2$, and $m - k_2\sqrt{n} \leq K \leq m - k_1\sqrt{n}$. Thus $R_n = P(K_F < m) + o(1)$.

**A.3. Proof of Lemma 2.** We want to show $\Delta(E_n, F_n, \mathcal{F}) \to 0$. We only show here that $\delta(F_n, E_n) \to 0$ as the other direction, namely, $\delta(E_n, F_n) \to 0$, is similar.

By the triangle inequality for deficiency we have

$$\delta(F_n, E_n) \leq \delta(F_n, F_{n+D\sqrt{n}}) + \delta(F_{n+D\sqrt{n}}, E_n).$$

From the assumption of Lemma 2 we know that the first term on the right-hand side goes to zero and hence it suffices to show that $\lim_{D\to\infty} \delta(F_{n+D\sqrt{n}}, E_n) \to 0$. Now let $\nu_n$ be a Poisson$(n+m)$ variable with $m = D\sqrt{n}$ and define $\nu_n^+ = \max(n, \nu_n)$. Let $F_{n+m}^\#$ be the experiment obtained by observing $x_1, x_2, \ldots, x_{\nu_n^+}$ i.i.d. with density $f$. Clearly $E_n \preceq F_{n+m}^\#$ (where $\preceq$ means "less informative" for experiments). We have

$$\Delta(F_{n+m}^\#, F_{n+m}, \mathcal{F}) \leq \|(\mathcal{L}(\nu_n^+), \mathcal{L}(\nu_n))\|_{\text{TV}} = P(\nu_n^+ \leq n - 1).$$

The Markov inequality gives

$$P(\nu_n^+ \leq n - 1) \leq \frac{\text{Var}(\nu_n^+)}{m^2} = \frac{(n+m)}{m^2} \leq \frac{2n}{m^2}.$$

This implies, since $m = D\sqrt{n}$, that

$$\delta(F_{n+D\sqrt{n}}, E_n) \leq \frac{2}{D^2}$$

and Lemma 2 follows on letting $D \to \infty$.

**A.4. Proof of Theorem 2.** As mentioned earlier, Theorem 2 is termed a globalization step in the asymptotic equivalence literature. The approach given here follows that of Nussbaum [13] and is by now somewhat standard. For simplicity of the notation we assume that $n$ is even. There are two steps.

*Step* 1. Split the observations $\{y_1, y_2, \ldots, y_n\}$ of $E_n$ into two sets of the same size,

$$\left\{\underline{y}_{1,n/2} = \left(y_{1i}; i = 1, \ldots, \frac{n}{2}\right), \underline{y}_{2,n/2} = \left(y_{2i}; i = 1, \ldots, \frac{n}{2}\right)\right\}.$$

Then define a new experiment $F_n^{\#}$ with the following independent observations:

$$\left\{\underline{y}_{1,n/2}, \left\{\underline{x}_{2,n/2}(\cdot) \text{ with intensity } \frac{n}{2}f\right\}\right\},$$

which is a modification of $E_n$ with the second set of observations in $E_n$ replaced by a set of observations from $F_n$. For ease of reading write $\underline{y}_1$, $\underline{y}_2$ and $\underline{x}_2$ to replace $\underline{y}_{1,n/2}$, $\underline{y}_{2,n/2}$ and $\underline{x}_{2,n/2}$. Let $\mathcal{F}_0 = \mathcal{F}(f_0, c_n)$ defined in (10). For any $\varepsilon > 0$, Lemma 3 tells us that the second set of observations in $E_n$ is locally asymptotically equivalent to a set of observations from $F_n$ uniformly in $f_0 \in \mathcal{F}$, that is, for all $f_0 \in \mathcal{F}$ there is a transition $K_{f_0}$ such that

$$\sup_{f_0 \in \mathcal{F}} \sup_{f \in \mathcal{F}_0} |K_{f_0} P_{2,f} - Q_{2,f}| \leq \varepsilon$$

when $n$ is sufficiently large, and from Proposition 9.2 in Nussbaum [13] every transition $K_{f_0}$ is given by a Markov kernel. With the first set of observations, we will construct an estimator $\hat{f}_n$ for $f$ that satisfies the optimality criterion given in (11). Because the parameter space is separable in Hellinger distance, and for fixed $n$ Hellinger distance is equivalent to the loss function $L_n$ by equation (12), we may further assume that $\hat{f}_n \in \mathcal{F}_0$ with $\mathcal{F}_0$ countable. Then one can show that $E_n$ and $F_n^{\#}$ are asymptotically equivalent. For any measurable set $B$ of experiment $F_n^{\#}$, define a randomization procedure

$$M(\underline{y}_1, \underline{y}_2, B) = \int 1_B((\underline{y}_1; \underline{x}_2)) K_{\hat{f}_n(\underline{y}_1)}(\underline{y}_2, d\underline{x}_2).$$

To show that $M$ is a Markov kernel, it is enough to check the measurability of $K_{\hat{f}_n(\underline{y}_1)}(\underline{y}_2, B_2)$ in $(\underline{y}_1, \underline{y}_2)$, for any given measurable set $B_2$. It is easy to see the measurability follows from the condition that $\hat{f}_n \in F_0^n$ with $F_0^n$ countable. Then

$$MP_f(B) = \int \int 1_B((\underline{y}_1; \underline{x}_2))(K_{\hat{f}_n(\underline{y}_1)} P_{2,f})(d\underline{x}_2) P_{1,f}(d\underline{y}_1),$$

which is expected to be close to

$$Q_f^{\#}(B) = \int \int 1_B((\underline{y}_1; \underline{x}_2)) Q_{2,f}(d\underline{x}_2) P_{1,f}(d\underline{y}_1).$$

Let $A_1 = \{\underline{y}_1, L_n(f, \hat{f}_{n,\varepsilon}) \leq \varepsilon\}$ and $\sup_f P(A_1^c) \leq \varepsilon$. Then

$$|M \cdot P_f(B) - Q_f^\#(B)|$$
$$= \left| \int \int 1_B((\underline{y}_1; \underline{x}_2))[(K_{\hat{f}_n(\underline{y}_1)} P_{2,f}) - Q_{2,f}](d\underline{x}_2) P_{1,f}(d\underline{y}_1) \right|$$
$$= \left| \int [1_{A_1^c}(\underline{y}_1) + 1_{A_1}(\underline{y}_1)] \right.$$
$$\left. \times \int 1_B((\underline{y}_1; \underline{x}_2))[(K_{\hat{f}_n(\underline{y}_1)} P_{2,f}) - Q_{2,f}](d\underline{x}_2) P_{1,f}(d\underline{y}_1) \right|$$
$$\leq 2P_{1,f}(A_1^c) + \sup_{f_0 \in \mathcal{F}} \sup_{f \in \mathcal{F}(f_0)} |K_{f_0} P_{2,f} - Q_{2,f}| \leq 3\varepsilon$$

uniformly over all $f \in \mathcal{F}$. Thus $\delta(E_n, F_n^\#) \leq 3\varepsilon$. Similarly, we can show $\delta(F_n^\#, E_n) \leq 3\varepsilon$. That is to say, $\Delta(E_n, F_n^\#, \mathcal{F}) \to 0$.

*Step* 2. We will then apply this procedure again to $F_n^\#$ in order to replace the first set of observations by its "asymptotically equivalent set," and obtain the compound experiment $F_n^{\#\#}$ where one observes two independent Poisson processes with intensity $\frac{n}{2}f$,

$$\{\underline{x}_{1,n}(\cdot), \underline{x}_{2,n}(\cdot)\}.$$

Here we need an estimator in $\mathcal{F}_0$ for $f$ which is derived from the second part of the observations in $F_n^\#$ and which has to satisfy the same optimality criterion. Similarly we have $\Delta(F_n^\#, F_n^{\#\#}, \mathcal{F}) \to 0$.

By applying a sufficiency argument, we see $F_n^{\#\#}$ is equivalent to $F_n$, so

$$\Delta(E_n, F_n, \mathcal{F}) \to 0.$$

Similarly we can show $\Delta(F_n, F_{n+D\sqrt{n}}, \mathcal{F}) \to 0$ using Lemma 3. Then the triangle inequality (18) for the deficiency distance gives

$$\Delta(E_n, E_{n+D\sqrt{n}}, \mathcal{F}) \leq \Delta(E_n, F_n, \mathcal{F}) + \Delta(F_n, F_{n+D\sqrt{n}}, \mathcal{F})$$
$$+ \Delta(F_{n+D\sqrt{n}}, E_{n+D\sqrt{n}}, \mathcal{F}) \to 0.$$

**A.5. Proof of Theorem 3.** The following simple lemma is used in the proof of Theorem 3.

LEMMA 4. *Let $N \sim \text{Poisson}(\lambda)$. Then*
$$\mathbb{P}\{|N - \lambda| \geq m_0\} \leq \exp(-m_0^3/(m_0 + \lambda)^2).$$

PROOF. The Chebyshev inequality gives

$$\mathbb{P}\{N - \lambda \geq m_0\} \leq \frac{E\exp(t(N-\lambda))}{\exp(tm_0)}$$
$$= \exp(\lambda(e^t - 1 - t) - tm_0) \qquad \text{for all } t \geq 0.$$

Let $t = m_0/(m_0 + \lambda)$. We know $e^t - 1 - t \leq t^2$ for $0 \leq t \leq 1$. Then
$$\mathbb{P}\{N - \lambda \geq m_0\} \leq \exp(\lambda m_0^2/(m_0 + \lambda)^2 - m_0^2/(m_0 + \lambda))$$
$$= \exp(-m_0^3/(m_0 + \lambda)^2).$$

Similarly we have $\mathbb{P}\{N - \lambda \leq -m_0\} \leq \exp(-m_0^3/(m_0 + \lambda)^2)$. □

PROOF OF THEOREM 3. For $\bar{f}_{(k)}$ defined in (6) and $c_n = \frac{1}{\sqrt{\log n}}$ set
$$B_n = \{x : |\bar{f}_{(k)}(x)| \leq 1/c_n, \ x \in [0,1]\}$$
and define $A_n$ by
$$A_n = \left\{\underline{x}_n : \sup_{x \in B_n} |\bar{f}_{(k)} - \widehat{\bar{f}}_n| \leq c_n\right\},$$
where $\hat{f}_n$ is the estimator defined in (14).

We first divide the expected loss into two pieces,
$$(19) \qquad E_f L_n(f, \hat{f}_n) = E_f 1_{A_n^c} L_n(f, \hat{f}_n) + E_f 1_{A_n} L_n(f, \hat{f}_n).$$

Note that the set $A_n$ can also be written as
$$A_n = \left\{\underline{x}_n : \max_{j \in \{j : J_{j,k}(f) \leq 1/c_n\}} \left|N_j - \frac{n}{k} J_{j,k}(f)\right| \leq c_n \frac{n}{k}\right\},$$
where $J_{j,k}$ is the averaging operator defined in Section 2.1. Since $N_j$ has a Poisson distribution it follows from Lemma 4 that
$$P\left(\left|N_j - \frac{n}{k} J_{j,k}(f)\right| \geq c_n \frac{n}{k}\right) \leq \exp\left(-\left(c_n \frac{n}{k}\right)^3 \bigg/ \left(c_n \frac{n}{k} + \frac{n}{k} J_{j,k}(f)\right)^2\right).$$

Since $c_n < \frac{1}{c_n}$ and $\frac{n}{k} = (\log n)^4$ it is easy to check that
$$\sup_{|J_{j,k}(f)| \leq 1/c_n} \left(c_n \frac{n}{k}\right)^3 \bigg/ \left(c_n \frac{n}{k} + \frac{n}{k} J_{j,k}(f)\right)^2 \geq \frac{1}{4} c_n^5 (\log n)^4 = \frac{1}{4}(\log n)^{3/2}.$$

Thus since $j$ only ranges from 1 to $k$ it immediately follows that
$$P(A_n^c) \leq k \exp(-\tfrac{1}{4}(\log n)^{3/2}) = o(n^{-\gamma}), \qquad \text{for any } \gamma > 1.$$

Now
$$L_n(f, \hat{f}_n) = \int_0^1 \frac{(\hat{f}_n - f)^2}{f + n^{-1/2}\hat{f}_n} \leq n^{1/2} \int_0^1 (\hat{f}_n + f) \leq n^{1/2}\left(\frac{N}{n} + 1\right).$$

Hence Cauchy–Schwarz yields
$$E_f 1_{A_n^c} L_n(f, \hat{f}_n) \leq n^{1/2} (P(A_n^c))^{1/2} \left(E\left(1 + \frac{N}{n}\right)^2\right)^{1/2} = o(1).$$

Now to bound the second term in (19) introduce sets $E_n$ and $F_n$ as follows:
$$E_n = \{x : |f - \bar{f}_{(k)}| \leq c_n/2, x \in [0,1]\}$$
and
$$F_n = \{x : \hat{f}_n = 0, x \in [0,1]\}.$$

Now the second term in (19) can be written as
$$E_f 1_{A_n} \int_{E_n^c \cup B_n^c} \frac{(f - \hat{f}_n)^2}{f + n^{-1/2}\hat{f}_n} + E_f 1_{A_n} \int_{E_n \cap B_n \cap F_n} \frac{(f - \hat{f}_n)^2}{f + n^{-1/2}\hat{f}_n}$$
$$+ E_f 1_{A_n} \int_{E_n \cap B_n \cap F_n^c} \frac{(f - \hat{f}_n)^2}{f + n^{-1/2}\hat{f}_n}$$
$$= R_1 + R_2 + R_3.$$

We take each of these terms one at a time. First it is convenient to break $R_1$ into two terms,
$$R_1 = E_f 1_{A_n} \int_{E_n^c} \frac{(f - \hat{f}_n)^2}{f + n^{-1/2}\hat{f}_n} + E_f 1_{A_n} \int_{B_n^c \cap E_n} \frac{(f - \hat{f}_n)^2}{f + n^{-1/2}\hat{f}_n} = R_{11} + R_{12}.$$

Now note that
$$R_{11} = E_f 1_{A_n} \int_{E_n^c \cap \{x : f(x) \leq 1/c_n\}} \frac{(f - \hat{f}_n)^2}{f + n^{-1/2}\hat{f}_n}$$
$$+ E_f 1_{A_n} \int_{E_n^c \cap \{x : f(x) > 1/c_n\}} \frac{(f - \hat{f}_n)^2}{f + n^{-1/2}\hat{f}_n}.$$

The definition of $\hat{f}_n$ with $0 \leq \hat{f}_n \leq 1/c_n$ then shows that
$$R_{11} \leq \sqrt{n} E_f \int_{E_n^c \cap \{x : f(x) \leq 1/c_n\}} (f + \hat{f}_n) + E_f \int_{E_n^c \cap \{x : f(x) > 1/c_n\}} \frac{f^2}{1/c_n}$$
$$\leq \frac{2\sqrt{n}}{c_n} \mu(E_n^c) + c_n \int f^2.$$

Now from $\int f^2 \leq C$ and the assumption (15) it follows that
$$R_{11} \leq C \left[\sqrt{n}\frac{2}{c_n}\right] k^{-\delta} + C c_n,$$
where $\delta > 1/2$ and $k = n/(\log n)^4$. Thus $R_{11} = o(1)$.

Now we consider $R_{12}$. It follows from the definition of $B_n$ and $E_n$ that on $B_n^c \cap E_n$ we have
$$f(x) \geq \bar{f}_{(k)}(x) - |f(x) - \bar{f}_{(k)}(x)| \geq 1/c_n - c_n/2 \geq \frac{1}{2c_n}.$$

Thus $R_{12}$ is bounded by

$$2c_n E_f \int (f - \hat{f}_n)^2 \leq 2c_n \int (f^2 + E_f \hat{f}_n^2) \leq 2c_n \int (f^2 + E_f \hat{\bar{f}}_n^2).$$

Simple moment calculations for a Poisson random variable give

$$\int E_f \hat{\bar{f}}_n^2 = \int \bar{f}_{(k)}^2 + \frac{k}{n} \leq \int f^2 + \frac{k}{n}.$$

Thus $R_{12} = o(1)$ since $\int f^2 \leq C$ by assumption. Hence $R_1 = R_{11} + R_{12} = o(1)$ uniformly over $\mathcal{F}$.

We now turn to $R_2$. Note that since $\hat{f}_n = 0$ on $F_n$,

$$R_2 = E_f 1_{A_n} \int_{E_n \cap B_n \cap F_n} \frac{(f - \hat{f}_n)^2}{f + n^{-1/2} \hat{f}_n} = E_f 1_{A_n} \int_{E_n \cap B_n \cap F_n} f.$$

Now note that since

$$f(x) \leq |f(x) - \bar{f}_{(k)}(x)| + |\bar{f}_{(k)}(x) - \hat{\bar{f}}_n(x)| + \hat{\bar{f}}_n(x),$$

it then follows from the definition of $E_n$, $B_n$ and $F_n$ that when $A_n$ occurs

$$f(x) \leq \frac{7c_n}{2}$$

and it immediately follows that $R_2 = o(1)$.

We finally turn to $R_3$. Since

$$f(x) \geq \hat{\bar{f}}_n(x) - |\hat{\bar{f}}_n(x) - \bar{f}_{(k)}(x)| - |\bar{f}_{(k)}(x) - f(x)|,$$

it then follows from the definition of $B_n$, $E_n$ and $F_n^c$ that when $A_n$ occurs $f(x) \geq c_n/2$ and since

$$|f(x) - \hat{\bar{f}}_n(x)| \leq |f(x) - \bar{f}_{(k)}(x)| + |\bar{f}_{(k)}(x) - \hat{\bar{f}}_n(x)|,$$

it follows that when $A_n$ occurs that and $x \in B_n \cap E_n$

$$|f - \hat{\bar{f}}_n| \leq 3c_n/2.$$

Hence $R_3$ is bounded by

$$\frac{2}{c_n} E_f 1_{A_n} \int_{E_n \cap B_n \cap F_n^c} (f - \hat{\bar{f}}_n)^2 = o(1).$$

The proof of (16) is complete since we have shown that $R_1 + R_2 + R_3 = o(1)$ uniformly over $\mathcal{F}$.

The proof of the theorem will be complete once we have shown that the assumptions of the theorem hold for Besov spaces $B_{p,q}^\alpha(M)$ with $\alpha p > 1/2$, $p \geq 1$ and $q \geq 1$. First note that if $\alpha - 1/p + 1/2 > 0$ and in particular if $\alpha p > 1/2$ and $p \geq 1$ then $B_{p,q}^\alpha(M)$ is compact in $L^2([0,1])$ (see the Appendix of Johnstone [6]) and so there is a constant $C$ such that $\int f^2 \leq C$ for all $f \in B_{p,q}^\alpha(M)$. Now the definition of the Besov norm shows that

$$2^{i\alpha}\|\bar{f}_{(2^{i+1})} - \bar{f}_{(2^i)}\|_p \leq M,$$

which implies for $p \geq 1$ and all $i_1 > i_0$

$$\|\bar{f}_{(2^{i_1+1})} - \bar{f}_{(2^{i_0})}\|_p \leq \sum_{i=i_0}^{i_1} \|\bar{f}_{(2^{i+1})} - \bar{f}_{(2^i)}\|_p \leq M \sum_{i=i_0}^{i_1} \left(\frac{1}{2^{i\alpha}}\right) \leq \frac{M}{(2^{i_0})^\alpha} \frac{2^\alpha}{2^\alpha - 1}.$$

Now take $k = 2^{i_0}$ and let $i_1 \to \infty$ to yield

$$k^{\alpha p} \int |f - \bar{f}_{(k)}|^p \leq \left(\frac{M 2^\alpha}{2^\alpha - 1}\right)^p.$$

Then the Chebyshev inequality gives

$$k^{\alpha p} c_k^p \mu\{x : |f - \bar{f}_{(k)}| > c_k\} \leq M_1^p.$$

Now let $c_k = \frac{1}{\sqrt{\log k}}$ to yield

$$\mu\left\{x : |f - \bar{f}_{(k)}| \frac{1}{\sqrt{\log k}}\right\} \leq M_1^p k^{-\alpha p} (\log k)^{p/2},$$

where $M_1 = M 2^\alpha / (2^\alpha - 1)$. Assumption (15) then clearly follows for $1/2 < \delta < \alpha p$. $\square$

DEPARTMENT OF STATISTICS
THE WHARTON SCHOOL
UNIVERSITY OF PENNSYLVANIA
PHILADELPHIA, PENNSYLVANIA 19104-6340
USA
E-MAIL: lowm@wharton.upenn.edu

DEPARTMENT OF STATISTICS
YALE UNIVERSITY
P.O.BOX 208290
NEW HAVEN, CONNECTICUT 06520-8290
USA
E-MAIL: huibin.zhou@yale.edu